\documentclass{amsart}

\usepackage{amsthm}
\usepackage{graphicx}
\usepackage{stmaryrd} 
\usepackage{tikz}
\usepackage{tikz-cd}
\usepackage{ bbold }
\usetikzlibrary{arrows}
\usepackage{verbatim}
\usepackage{amssymb,amsfonts}
\usepackage[all,arc]{xy}
\usepackage{enumerate}
\usepackage{mathrsfs, mathabx}
\usepackage{xcolor}[2007/01/21]
\usepackage{hyperref}
\usepackage{soul}
\usepackage{esint} 
\usepackage[all]{xy}
\usepackage[margin=1.0in]{geometry}
\usepackage[
   open,
   openlevel=2,
   atend
 ]{bookmark}[2011/12/02]
\usepackage{graphics}

\bookmarksetup{color=blue}

\theoremstyle{definition} 
\newtheorem{thm}{Theorem}[section]

\newtheorem{lem}[thm]{Lemma}

\newtheorem{quest}[thm]{Question}

\theoremstyle{definition}

\newtheorem{notn}[thm]{Notation}

\newtheorem{rmk}[thm]{Remark}

\theoremstyle{definition}

\theoremstyle{remark}

\newcommand{\GL}{\mathrm{GL}}

\newcommand{\Aut}{\mathrm{Aut}}

\newcommand{\bZ}{\mathbb{Z}}
\newcommand{\bQ}{\mathbb{Q}}
\newcommand{\bR}{\mathbb{R}}

\newcommand{\bA}{\mathbb{A}}

\newcommand{\bF}{\mathbb{F}}

\newcommand{\mc}{\mathcal}
\newcommand{\sm}{\smallsetminus}

\newcommand{\Mat}{\mathrm{Mat}}
\newcommand{\Prob}{\mathrm{Prob}}

\newcommand{\ld}{\lambda}
\newcommand{\al}{\alpha}

\newcommand{\ol}{\overline}

\newcommand{\ra}{\rightarrow}

\newcommand{\hra}{\hookrightarrow}

\newcommand{\lt}{\left}
\newcommand{\rt}{\right}

\newcommand{\op}{\oplus}
\newcommand{\acts}{\lefttorightarrow}

\newcommand{\be}{\begin{enumerate}}
\newcommand{\ee}{\end{enumerate}}

\newcommand{\bi}{\begin{itemize}}
\newcommand{\ei}{\end{itemize}}

\newcommand{\bbm}{\begin{bmatrix}}
\newcommand{\ebm}{\end{bmatrix}}

\numberwithin{equation}{section}

\begin{document}

\title{Jordan--Landau theorem for matrices over finite fields}
\date{\today}
\author{Gilyoung Cheong, Jungin Lee, Hayan Nam, and Myungjun Yu}
\address{G. Cheong -- Department of Mathematics, University of California, Irvine, 340 Rowland Hall, Irvine, CA 92697 \newline
J. Lee -- School of Mathematics, Korea Institute for Advanced Study, Seoul 02455, Korea \newline
H. Nam -- Department of Mathematics, Duksung Women's University, 33 Samyang-ro 144-gil, Seoul, South Korea \newline
M. Yu -- Department of Mathematics, Yonsei University, Seoul 03722, South Korea}
\email{gilyounc@uci.edu, jilee@kias.re.kr, hnam@duksung.ac.kr, mjyu@yonsei.ac.kr}

\begin{abstract}
Given a positive integer $r$ and a prime power $q$, we estimate the probability that the characteristic polynomial $f_{A}(t)$ of a random matrix $A$ in $\GL_{n}(\bF_{q})$ is square-free with $r$ (monic) irreducible factors when $n$ is large. We also estimate the analogous probability that $f_{A}(t)$ has $r$ irreducible factors counting with multiplicity. In either case, the main term $(\log n)^{r-1}((r-1)!n)^{-1}$ and the error term $O((\log n)^{r-2}n^{-1})$, whose implied constant only depends on $r$ but not on $q$ nor $n$, coincide with the probability that a random permutation on $n$ letters is a product of $r$ disjoint cycles. The main ingredient of our proof is a recursion argument due to S. D. Cohen, which was previously used to estimate the probability that a random degree $n$ monic polynomial in $\bF_{q}[t]$ is square-free with $r$ irreducible factors and the analogous probability that the polynomial has $r$ irreducible factors counting with multiplicity. We obtain our result by carefully modifying Cohen's recursion argument in the matrix setting, using Reiner's theorem that counts the number of $n \times n$ matrices with a fixed characteristic polynomial over $\bF_{q}$.
\end{abstract}

\maketitle


\section{Introduction}

\hspace{3mm} It is folklore that the cycle decomposition of a random permutation is analogous to the prime decomposition of a random positive integer. For instance, Granville \cite[Section 2.1]{Gra} notes that given $r \in \bZ_{\geq 1}$, a theorem of Jordan \cite[p.161]{Jor} says
\[\underset{\sigma \in S_{n}}{\Prob}\left(\begin{array}{c}
\sigma \text{ is a product of}\\ 
r \text{ disjoint cycles}
\end{array}\right) \sim \frac{(\log n)^{r-1}}{(r-1)! n}\]

for large $n \in \bZ_{\geq 1}$, while a theorem of Landau \cite[p.211]{Lan}, which can also be found in \cite[p.491, Theorem 437]{HW}, says
\begin{align*}
\underset{1 \leq N \leq x}{\Prob}\left(\begin{array}{c}
N \text{ is square-free with} \\
r \text{ prime factors}
\end{array}\right) &\sim \underset{1 \leq N \leq x}{\Prob}\left(\begin{array}{c}
N \text{ has } r \text{ prime factors} \\
\text{counting with multiplicity}
\end{array}\right)\\ 
&\sim \frac{(\log \log x)^{r-1}}{(r-1)!\log x}
\end{align*}

for large $x \in \mathbb{R}_{\geq 1}$.\footnote{In the literature, the number of $\sigma \in S_{n}$ with $r$ disjoint cycles is called a \textbf{Stirling number of the first kind}. We write ``$f(x) \sim g(x)$ for large $x$'' to mean that $f(x)/g(x) \ra 1$ as $x \ra \infty$.} The first probability is given uniformly at random from the set $S_{n}$ of permutations on $n$ letters. The second probability is given by choosing the integer $N$ uniformly at random from the set $\{1, 2, \dots, \lfloor x \rfloor\}$. Landau's theorem is a generalization of the Prime Number Theorem (i.e., the case $r = 1$). The function field analogue of Laudau's theorem is that given $r \in \bZ_{\geq 1}$, for large $n$, with a fixed prime power $q$, we have
\begin{align*}
\underset{f \in \bA^{n}(\bF_{q})}{\Prob}
\left(\begin{array}{c}
f(t) \text{ is square-free with} \\
r\text{ irreducible factors}
\end{array}\right) &\sim \underset{f \in \bA^{n}(\bF_{q})}{\Prob}
\left(\begin{array}{c}
f(t) \text{ has } r \text{ irreducible factors} \\
\text{counting with multiplicity}
\end{array}\right) \\
&\sim \frac{(\log n)^{r-1}}{(r-1)! n},
\end{align*}

where we denote by $\bA^{n}(\bF_{q})$ the set of degree $n$ monic polynomials in $\bF_{q}[t]$ over the finite field $\bF_{q}$ of $q$ elements.\footnote{Throughout the paper, we write ``irreducible factors'' to mean monic irreducible factors.} We note that the estimate, which is due to Cohen \cite[Theorem 6]{Coh}, does not depend on $q$. Car \cite{Car} also proved the above estimate with a different method, and both papers of Cohen and Car prove that the desired probabilities are equal to
\[\frac{(\log n)^{r-1}}{(r-1)! n} + O\left(\frac{(\log n)^{r-2}}{n}\right),\]

where the implied constant only depends on $r$ but not on $q$ nor $n$. In this paper, we give a matrix version of this theorem. We asymptotically compute, given a (uniformly) random matrix $A \in \GL_{n}(\bF_{q})$, the probability that the characteristic polynomial $f_{A}(t) = \det(tI_{n} - A)$ is square-free with $r$ irreducible factors. We also estimate the analogous probability that $f_{A}(t)$ has $r$ irreducible factors counting with multiplicity. The same method gives us analogous results when we replace $\GL_{n}(\bF_{q})$ with $\Mat_{n}(\bF_{q})$, the set of all $n \times n$ matrices over $\bF_{q}$ including singular matrices.

\begin{thm}\label{main} Let $n, r \in \bZ_{\geq 1}$, and fix a prime power $q$. We have
\begin{align*}
&\underset{A \in \GL_{n}(\bF_{q})}{\Prob}
\left(\begin{array}{c}
f_{A}(t) \text{ is square-free with} \\
r\text{ irreducible factors}
\end{array}\right) = \frac{(\log n)^{r-1}}{n(r-1)!} + O\left( \frac{(\log n)^{r-2}}{n} \right)
\end{align*}

and
\begin{align*}
&\underset{A \in \GL_{n}(\bF_{q})}{\Prob}
\left(\begin{array}{c}
f_{A}(t) \text{ has } r \text{ irreducible factors}\\
\text{counting with multiplicity}
\end{array}\right) = \frac{(\log n)^{r-1}}{n(r-1)!} + O\left( \frac{(\log n)^{r-2}}{n} \right),
\end{align*}

where the implied constants may depend on $r$ but not on $q$ nor $n$. We also have
\begin{align*}
&\underset{A \in \Mat_{n}(\bF_{q})}{\Prob}
\left(\begin{array}{c}
f_{A}(t) \text{ is square-free with} \\
r\text{ irreducible factors}
\end{array}\right) = \left(\frac{(\log n)^{r-1}}{n(r-1)!} + O\left( \frac{(\log n)^{r-2}}{n} \right)\right)\prod_{i=1}^{n}(1 - q^{-i})
\end{align*}

and
\begin{align*}
&\underset{A \in \Mat_{n}(\bF_{q})}{\Prob}
\left(\begin{array}{c}
f_{A}(t) \text{ has } r \text{ irreducible factors} \\
\text{counting with multiplicity}
\end{array}\right) = \left(\frac{(\log n)^{r-1}}{n(r-1)!} + O\left( \frac{(\log n)^{r-2}}{n} \right)\right)\prod_{i=1}^{n}(1 - q^{-i}),
\end{align*}

where the implied constants may depend on $r$ but not on $q$ nor $n$.
\end{thm}

\begin{rmk} Theorem \ref{main} is our main result of this paper. Our contributions are the estimates for fixed $q$. We find it surprising that the implied constants for the big-O terms in the theorem do not depend on $q$ nor $n$. Of course, the same big-O term also appears in Cohen's result for random monic polynomials, and our proof makes use of a modified version of Cohen's argument. However, we still do not have a satisfying philosophical understanding on why this independence on $q$ and $n$ happens other than just knowing that the computations work out.
\end{rmk}

\begin{rmk}
The factor $\prod_{i=1}^{n}(1 - q^{-i})$ in statement for $\Mat_{n}(\bF_{q})$ is the fraction $|\GL_{n}(\bF_{q})|/|\Mat_{n}(\bF_{q})|$, which happens to occur when we work with all matrices instead of just invertible ones. Roughly speaking, such a factor arises because the characteristic polynomial is invariant under the conjugate action of $\GL_{n}(\bF_{q})$, so various formulas in our arguments end up having $|\GL_{n}(\bF_{q})|$, which is canceled when we consider probabilities with the sample space $\GL_{n}(\bF_{q})$ instead of $\Mat_{n}(\bF_{q})$.
\end{rmk}

\

\subsection{Related work and further questions} Perhaps, it is possible to generalize Theorem \ref{main} further. An interesting clue can be found in a result of Hansen and Schmutz \cite[Theorem 1]{HS}, who used the cycle index of $\GL_{n}(\bF_{q})$ to show that there exist $c_{1}, c_{2} \in \bR_{> 0}$ such that for any $m \in \bZ_{\geq 1}$ satisfying $c_{1} \log(n) \leq m \leq n$ and for any choice of $r_{m+1}, r_{m+2}, \dots, r_{n} \in \bZ_{\geq 0}$, we have
\[\lt| \underset{A \in \GL_{n}(\bF_{q})}{\Prob}
\left(\begin{array}{c}
f_{A}(t) \text{ has } r_{j} \text{ degree } j \\
\text{irreducible factors for } m+1 \leq j \leq n \\
\text{counting with multiplicity}
\end{array}\right) - \underset{f \in \bA^{n}(\bF_{q})}{\Prob}
\left(\begin{array}{c}
f(t) \text{ has } r_{j} \text{ degree } j \\
\text{irreducible factors for } m+1 \leq j \leq n \\
\text{counting with multiplicity}
\end{array}\right) \rt| \leq \frac{c_{2}}{m}.\]

\hspace{3mm} This brings us a natural question:

\begin{quest}\label{Q} Given any distinct $d_{1}, \dots, d_{l} \in \bZ_{\geq 1}$ and not necessarily distinct $r_{1}, \dots, r_{l} \in \bZ_{\geq 0}$, do we have a reasonably good estimate for
\[\lt| \underset{A \in \GL_{n}(\bF_{q})}{\Prob}
\left(\begin{array}{c}
f_{A}(t) \text{ has } r_{j} \text{ degree } d_{j} \\
\text{irreducible factors for } 1 \leq j \leq l \\
\text{counting with multiplicity}
\end{array}\right) 
- \underset{f \in \bA^{n}(\bF_{q})}{\Prob}
\left(\begin{array}{c}
f(t) \text{ has } r_{j} \text{ degree } d_{j} \\
\text{irreducible factors for } 1 \leq j \leq l \\
\text{counting with multiplicity}
\end{array}\right) \rt|?\]

What about
\[\lt| \underset{A \in \GL_{n}(\bF_{q})}{\Prob}
\left(\begin{array}{c}
f_{A}(t) \text{ has } r_{j} \text{ degree } d_{j} \\
\text{irreducible factors for } 1 \leq j \leq l \\
\text{counting with multiplicity}
\end{array}\right) 
- \underset{\sigma \in S_{n}}{\Prob}\left(\begin{array}{c}
\sigma \text{ has } r_{j} \text{ disjoint} \\ 
d_{j}\text{-cycles for } 1 \leq j \leq l
\end{array}\right) \rt|?\]
\end{quest}

\begin{rmk} It turns out that as $q \ra \infty$ both quantities in Question \ref{Q} go to $0$, and this is well-known among experts. Hence, it is interesting consider this question with fixed $q$. (More details on the large $q$ limit can be found in an older version of this paper \cite{CNY}.\footnote{We decided to omit this proof for the current version of our paper because it is only relevant to Question \ref{Q}.}) By ``a reasonably good estimate'', we mean a sequence $E_{n} \in \bR_{\geq 0}$ that only depends on $n$ such that $\lim_{n \ra \infty}E_{n} = 0$ and
\[\lt| \underset{A \in \GL_{n}(\bF_{q})}{\Prob}
\left(\begin{array}{c}
f_{A}(t) \text{ has } r_{j} \text{ degree } d_{j} \\
\text{irreducible factors for } 1 \leq j \leq l \\
\text{counting with multiplicity}
\end{array}\right) 
- \underset{f \in \bA^{n}(\bF_{q})}{\Prob}
\left(\begin{array}{c}
f(t) \text{ has } r_{j} \text{ degree } d_{j} \\
\text{irreducible factors for } 1 \leq j \leq l \\
\text{counting with multiplicity}
\end{array}\right) \rt| \leq E_{n}.\]

The result of Hansen and Schmutz provides a partial answer to Question \ref{Q}, but it is still far from answering the full question, which is an interesting starting point for a future work.
\end{rmk}

\

\subsection{Outline of the proof of Theorem \ref{main}} For the square-free case, we prove Theorem \ref{main} by carefully modifying Cohen's arguments in \cite{Coh} that proves 
\[\underset{f \in \bA^{n}(\bF_{q})}{\Prob}
\left(\begin{array}{c}
f(t) \text{ is square-free with} \\
r\text{ irreducible factors}
\end{array}\right) = \frac{(\log n)^{r-1}}{(r-1)!n} + O\left(\frac{(\log n)^{r-2}}{n}\right),\]

where the implied constant only depends on $r$ but not on $q$ nor $n$. More specifically, we first derive a recursive formula for
\[\#\left\{\begin{array}{c}
A \in \Mat_{n}(\bF_{q}) : 
f_{A}(t) \text{ is square-free with} \\
r\text{ irreducible factors}
\end{array}
\right\}\]

in $r \geq 1$ and then estimate the count inductively, using the base case $r = 1$ given as follows
\begin{align*}
\#\left\{\begin{array}{c}
A \in \Mat_{n}(\bF_{q}) : 
f_{A}(t) \text{ is square-free with} \\
1\text{ irreducible factor}
\end{array}
\right\} &= \#\{
A \in \Mat_{n}(\bF_{q}) : 
f_{A}(t) \text{ is irreducible}\} \\
&= \frac{|\GL_{n}(\bF_{q})|M_{q}(n)}{q^{n}-1},
\end{align*}

where $M_{q}(n)$ is the number of monic irreducible polynomials of degree $n$ in $\bF_{q}[t]$. We recall that
\[M_{q}(n) = \frac{1}{n} \sum_{d | n} \mu(d)q^{n/d}\]

for any $n \geq 1$, where $\mu(d)$ is the M\"obius function. From this formula, one can deduce that (e.g., see \cite[Lemma 4]{Pol})
\[\frac{q^{n}}{n} - \frac{2q^{n/2}}{n}\leq M_{q}(n) \leq \frac{q^{n}}{n},\]

which we shall use in the rest of our paper without mentioning again. For the non-square-free case, we use a similar strategy, but deriving a recursive formula in $r \geq 1$ for  
\[\#\left\{\begin{array}{c}
A \in \Mat_{n}(\bF_{q}) : 
f_{A}(t) \text{ has } r \text{ irreducible factors} \\
\text{counting with multiplicity}
\end{array}
\right\}\]

gets more difficult because given any monic polynomial $g(t) \in \bA^{n}(\bF_{q})$, the formula for the number of $A \in \Mat_{n}(\bF_{q})$ such that $f_{A}(t) = g(t)$ is more complicated when $g(t)$ is not square-free. Nevertheless, we manage to get a recursive formula for the non square-free case, and similar but more convoluted arguments to the square-free case yield the result for the non-square-free case.

\

\section{Number of matrices with a fixed characteristic polynomial}

\hspace{3mm} In this section, we fix a prime power $q$.  In order to derive recursive formulas we need for Theorem \ref{main}, it is crucial for us to know the number of $n \times n$ matrices in $\Mat_{n}(\bF_{q})$ with a fixed characteristic polynomial for any given $n \in \bZ_{\geq 1}$. This number is explicitly known due to Reiner \cite{Rei}. We denote by $\mc{P}_{q}$ the set of all monic irreducible polynomials in $\bF_{q}[t]$.

\begin{thm}[Reiner]\label{matchar} Let $n, l \in \bZ_{\geq 1}$ and consider any distinct $P_{1}, \dots, P_{l} \in \mc{P}_{q}$ and not necessarily distinct $n_{1}, \dots, n_{l} \in \bZ_{\geq 1}$ such that $n_{1}\deg(P_{1}) + \cdots + n_{l}\deg(P_{l}) = n$. Then
\[\#\{A \in \Mat_{n}(\bF_{q}) : f_{A}(t) = P_{1}(t)^{n_{1}}P_{2}(t)^{n_{2}} \cdots P_{l}(t)^{n_{l}}\} = \frac{|\GL_{n}(\bF_{q})|(q^{\deg(P_{1})})^{n_{1}^{2} - n_{1}}(q^{\deg(P_{2})})^{n_{2}^{2} - n_{2}} \cdots (q^{\deg(P_{l})})^{n_{l}^{2} - n_{l}}}{|\GL_{n_{1}}(\bF_{q^{\deg(P_{1})}})| |\GL_{n_{2}}(\bF_{q^{\deg(P_{2})}})| \cdots |\GL_{n_{l}}(\bF_{q^{\deg(P_{l})}})|}.\]
\end{thm}

\subsection{Two proofs of Theorem \ref{matchar}} We include two different proofs of Theorem \ref{matchar} we found in the literature to help the reader compare the proofs. Nevertheless, any reader who is not interested in the proof of Theorem \ref{matchar} is welcome to skip this subsection.

\hspace{3mm} For any $A \in \Mat_{n}(\bF_{q})$ such that $f_{A}(t) = P_{1}(t)^{n_{1}}P_{2}(t)^{n_{2}} \cdots P_{l}(t)^{n_{l}}$, the $\bF_{q}[t]$-module structure on $\bF_{q}^{n}$, given by setting the multiplication by $t$ as the left-multiplication of $A$, can be decomposed as
\[(A \acts \bF_{q}^{n}) = A[P_{1}^{n_{1}}] \op \cdots \op A[P_{l}^{n_{l}}],\]

where
\[A[P_{i}^{n_{i}}] := \{v \in \bF_{q}^{n} : P_{i}(A)^{n_{i}}v = 0\},\]

which is a module over $\bF_{q}[t]/(P_{i}(t)^{n_{i}})$. If $V_{i}$ is a module over $\bF_{q}[t]/(P_{i}(t)^{n_{i}})$ with $\bF_{q}$-dimension $n_{i}\deg(P_{i})$ for $1 \leq i \leq l$, then applying the orbit-stabilizer theorem on the conjugate action of $\GL_{n}(\bF_{q})$ on $\Mat_{n}(\bF_{q})$, we see that the number of $A \in \Mat_{n}(\bF_{q})$ such that $(A \acts \bF_{q}^{n}) \simeq V_{1} \op \cdots \op V_{l}$ is 
\[\frac{|\GL_{n}(\bF_{q})|}{|\Aut_{\bF_{q}[t]}(V_{1})| |\Aut_{\bF_{q}[t]}(V_{2})| \cdots |\Aut_{\bF_{q}[t]}(V_{l})|},\]

where $\Aut_{\bF_{q}[t]}(V_{i})$ is the group of $\bF_{q}[t]$-automorphisms of $V_{i}$. Thus, we have
\begin{align*}
\#\{A \in \Mat_{n}(\bF_{q}) : f_{A}(t) = P_{1}(t)^{n_{1}}P_{2}(t)^{n_{2}} \cdots P_{l}(t)^{n_{l}}\} &= \sum_{[V_{1}], \dots, [V_{l}]} \frac{|\GL_{n}(\bF_{q})|}{|\Aut_{\bF_{q}[t]}(V_{1})| |\Aut_{\bF_{q}[t]}(V_{2})| \cdots |\Aut_{\bF_{q}[t]}(V_{l})|} \\
&= |\GL_{n}(\bF_{q})|\prod_{i=1}^{l}\sum_{\substack{[V_{i}] : \\ \dim(V_{i}) = n_{i}\deg(P_{i})}} \frac{1}{|\Aut_{\bF_{q}[t]}(V_{i})|},
\end{align*}

where $[V_{i}]$ varies among the isomorphism classes of $\bF_{q}[t]/(P_{i}(t))^{n_{i}})$-modules such that $\dim_{\bF_{q}}(V_{i}) = n_{i}\deg(P_{i})$. When $n_{1} = \cdots = n_{l} = 1$, which is the only case we need for the square-free case of Theorem \ref{main}, the result trivially follows from here, but the desired formula is more difficult to derive when any $n_{i} \geq 2$.

\hspace{3mm} In any case, Theorem \ref{matchar} reduces to
\begin{equation}\label{Rei}
\sum_{\substack{[V_{i}] : \\ \dim(V_{i}) = n_{i}\deg(P_{i})}} \frac{1}{|\Aut_{\bF_{q}[t]}(V_{i})|} = \frac{(q^{\deg(P_{i})})^{n_{i}^{2} - n_{i}}}{|\GL_{n_{i}}(\bF_{q^{\deg(P_{i})}})|},
\end{equation}

and Reiner \cite{Rei} shows that a combinatorial statement about a sum over certain partitions due to Fine and Herstein \cite{FH} directly implies (\ref{Rei}), which completes his proof of Theorem \ref{Rei}.

\hspace{3mm} Gerstenhaber \cite{Ger} gives a more algebraic proof of (\ref{Rei}), which we explain for the rest of this subsection. Let $d_{i} := \deg(P_{i})$. Choose a root $\ld \in \ol{\bF_{q}}$ of $P_{i}(t) \in \bF_{q}[t]$, which gives rise to an embedding
\[\iota : \bF_{q^{d_i}} = \bF_{q}[\ld] \hra \Mat_{d_i}(\bF_{q})\] 

of (generally non-commutative) $\bF_{q}$-algebras by mapping $\ld$ to the companion matrix, say $M$, of $P_{i}(t)$ over $\bF_{q}$. Consider the induced embedding $\tilde{\iota} : \Mat_{n_{i}}(\bF_{q^{d_{i}}}) \hra \Mat_{n_{i}}(\Mat_{d_{i}}(\bF_{q}))$ and observe that every matrix $X \in \Mat_{d_{i}n_{i}}(\bF_{q}) = \Mat_{n_{i}}(\Mat_{d_{i}}(\bF_{q}))$ such that $P_{i}(X)^{n_{i}} = 0$ must be similar to a matrix of the form $\tilde{\iota}(\ld I_{n_{i}} + E)$, where $E \in \Mat_{n_{i}}(\bF_{q^{d_{i}}})$ is a nilpotent matrix. Given any $C \in \GL_{d_{i}m_{i}}(\bF_{q})$ such that $C^{-1}\tilde{\iota}(\ld I_{n_{i}} + E)C = \tilde{\iota}(\ld I_{n_{i}} + E')$ for some nilpotent matrices $E, E' \in \Mat_{n_{i}}(\bF_{q^{d_{i}}})$, we have
\begin{equation}\label{comm}
C^{-1}\tilde{\iota}(\ld I_{n_{i}})C = C^{-1}\tilde{\iota}(\ld I_{n_{i}} + E)^{q^{d_{i}n_{i}}}C = \tilde{\iota}(\ld I_{n_{i}} + E')^{q^{d_{i}n_{i}}} = \tilde{\iota}(\ld I_{n_{i}}).
\end{equation}

Note that any $B \in \Mat_{d_{i}}(\bF_{q})$ is in the image of $\iota$ if and only if $BM = MB$. This is because the characteristic polynomial of $M$ is $P(t) \in \bF_{q}[t]$ which is irreducible. Indeed, this implies that the $\bF_{q}[t]$-module structure given by the left-multiplication $M \acts \bF_{q}^{n}$ is cyclic while $B$ is an endomorphism for such a module. A similar reasoning shows that any $\tilde{B} \in \Mat_{d_{i}n_{i}}(\bF_{q}) = \Mat_{n_{i}}(\Mat_{d_{i}}(\bF_{q}))$ is in the image of $\tilde{\iota}$ if and only if $\tilde{B}\tilde{\iota}(\ld I_{d_{i}}) = \tilde{\iota}(\ld I_{d_{i}})\tilde{B}$. Thus, we see (\ref{comm}) implies that $C$ is in the image of $\tilde{\iota} : \Mat_{n_{i}}(\bF_{q^{d_{i}}}) \hra \Mat_{n_{i}}(\Mat_{d_{i}}(\bF_{q}))$. From the definition of $\tilde{\iota}$, it follows that any $C \in \GL_{d_{i}n_{i}}(\bF_{q})$ such that $C^{-1}\tilde{\iota}(\ld I_{n_{i}} + E)C = \tilde{\iota}(\ld I_{n_{i}} + E')$ for some nilpotent matrices $E, E' \in \Mat_{n_{i}}(\bF_{q^{d_{i}}})$ must be mapped from $\GL_{n_{i}}(\bF_{q^{d_{i}}})$ under $\tilde{\iota}$. The number of nilpotent matrices in $\Mat_{n_{i}}(\bF_{q^{d_{i}}})$ is $(q^{d_{i}})^{n_{i}^{2} - n_{i}}$ by Fine and Herstein \cite{FH}, so the number of $X \in \Mat_{d_{i}n_{i}}(\bF_{q})$ such that $P_{i}(X)^{n_{i}} = 0$ (which necessarily implies that $f_{X}(t) = P_{i}(X)^{n_{i}}$) is $(q^{d_{i}})^{n_{i}^{2} - n_{i}}|\GL_{d_{i}n_{i}}(\bF_{q})|/|\GL_{n_{i}}(\bF_{q^{d_{i}}})|$. Applying the orbit-stabilizer theorem to the conjugate action of $\GL_{d_{i}n_{i}}(\bF_{q})$ on the set of such $X \in \Mat_{d_{i}n_{i}}(\bF_{q})$, this implies that
\[\sum_{\substack{[V_{i}] : \\ \dim(V_{i}) = n_{i}d_{i}}} \frac{|\GL_{d_{i}n_{i}}(\bF_{q})|}{|\Aut_{\bF_{q}[t]}(V_{i})|} = \frac{(q^{d_{i}})^{n_{i}^{2} - n_{i}}|\GL_{d_{i}n_{i}}(\bF_{q})|}{|\GL_{n_{i}}(\bF_{q^{d_{i}}})|},\]

which implies (\ref{Rei}).

\

\section{Recursion}

\hspace{3mm} In this section, we fix a prime power $q$. 

\subsection{Square-free case} In this subsection, we provide a recursive formula for
\[\xi_{r}(n) := \#\left\{\begin{array}{c}
A \in \Mat_{n}(\bF_{q}) : 
f_{A}(t) \text{ is square-free with} \\
r\text{ irreducible factors}
\end{array}
\right\}.\]

in $r \in \bZ_{\geq 1}$ for fixed $n \in \bZ_{\geq 0}$, by mimicking the argument given for \cite[Theorem 4]{Coh}. Note that $\GL_{0}(\bF_{q}) = \Mat_{0}(\bF_{q})$ is the set of $\bF_{q}$-linear maps of the trivial vector space over $\bF_{q}$, so $\xi_{r}(0) = 0$ for any $r \in \bZ_{\geq 1}$. It also follows that
\[\xi_{0}(n) = \left\{
	\begin{array}{ll}
	0 \mbox{ if } n \geq 1,  \\
	1 \mbox{ if } n = 0,
	\end{array}\right.\]

For convenience, we define
\[\xi_{r}(n) := 0\]

whenever $r < 0$ or $n < 0$. Note that in any case, we have $\xi_{r}(n) = 0$ if $r > n$.

\hspace{3mm} Recall that for any $n \in \bZ_{\geq 1}$, we have
\[|\GL_{n}(\bF_{q})| = (q^{n} - 1)(q^{n} - q) \cdots (q^{n} - q^{n-1}) = \prod_{i=0}^{n-1}(q^{n} - q^{i}).\]

When $n = 0$, we have $|\GL_{0}(\bF_{q})| = 1$. For our purpose, we would also like to have the version of $|\GL_{n}(\bF_{q})|$ for negative $n$, which does not make sense in linear algebra because there are no negative dimensional vector spaces. However, the product $\prod_{i=0}^{n-1}(q^{n} - q^{i})$ makes sense for negative $n$, which would give us the empty product, resulting in $1$. Thus, we define
\[|\GL_{n}(\bF_{q})| := 1\]

for $n \in \bZ_{< 0}$.

\hspace{3mm} We now state our recursion formula for $\xi_{r}(n)$. Since $q$ is fixed, we write $\pi(n) := M_{q}(n)$ to mean the number of monic irreducible polynomials of degree $n$ in $\bF_{q}[t]$ for ease of notation.

\begin{thm}\label{recursion} Given $n, r \in \bZ_{\geq 0}$, we have
\[n\xi_{r}(n) =  |\GL_{n}(\bF_{q})|\sum_{s=1}^{n}\sum_{e=1}^{r}(-1)^{e-1} \frac{s\pi(s)\xi_{r-e}(n-es)}{(q^{s}-1)^{e} |\GL_{n - es}(\bF_{q})|}.\]

Equivalently, writing $\phi_{r}(n) := \xi_{r}(n)/|\GL_{n}(\bF_{q})|$, we have
\[n\phi_{r}(n) =  \sum_{s=1}^{n}\sum_{e=1}^{r}(-1)^{e-1} \frac{s\pi(s)\phi_{r-e}(n-es)}{(q^{s}-1)^{e}}.\]
\end{thm}

\begin{rmk} We note the resemblance of our formula for $\phi_{r}(n)$ given in Theorem \ref{recursion} and the one for
\[\pi_{r}(n) := \#\left\{\begin{array}{c}
f \in \bA^{n}(\bF_{q}) : 
f(t) \text{ is square-free with} \\
r\text{ irreducible factors}
\end{array}
\right\}\]

given in \cite[Theorem 4]{Coh}:
\[n\pi_{r}(n) = \sum_{s=1}^{n}\sum_{e=1}^{r}(-1)^{e-1}s \pi(s) \pi_{r-e}(n-es).\]
\end{rmk}

\hspace{3mm} Recall that we denote by $\mc{P}_{q}$ the set of all monic irreducible polynomials in $\bF_{q}[t]$. Given $P \in \mc{P}_{q}$ and $r, n \in \bZ_{\geq 0}$, we write 
\[\psi_{r}(n ; P) := \#\left\{\begin{array}{c}
A \in \Mat_{n}(\bF_{q}) : 
f_{A}(t) \text{ is square-free with} \\
r\text{ irreducible factors and } P(t) \nmid f_{A}(t)
\end{array}
\right\}\]

Note that
\bi
	\item if $n \geq 1$, we have $\psi_{0}(n ; P) = 0$;
	\item if $r \geq 1$, we have $\psi_{r}(0 ; P) = 0$; 
	\item $\psi_{0}(0 ; P) = 1$. 
\ei

We also define 
\bi
	\item $\psi_{r}(n ; P) := 0$ whenever $r < 0$ or $n < 0$
\ei

for convenience. In any case, we have $\psi_{r}(n, P) = 0$ if $r < n$. For $n \in \bZ_{\geq 1}$, we note that 
\[\xi_{r}(n) - \psi_{r}(n ; P) = \#\left\{\begin{array}{c}
A \in \Mat_{n}(\bF_{q}) : 
f_{A}(t) \text{ is square-free with} \\
r\text{ irreducible factors and } P(t) | f_{A}(t)
\end{array}
\right\}\]

The following lemma contains most of the proof of Theorem \ref{recursion}:

\begin{lem}\label{rec} Let $P \in \mc{P}_{q}$ with $\deg(P) = s$. For any $n, r  \in \bZ_{\geq 0}$,  we have
\[\xi_{r}(n) - \psi_{r}(n; P) = |\GL_{n}(\bF_{q})|\sum_{e=1}^{r}(-1)^{e-1}\frac{\xi_{r-e}(n-es)}{|\GL_{n-es}(\bF_{q})|(q^{s}-1)^{e}}.\]
\end{lem}

\begin{proof} When $r = 0$, both sides are $0$, so we only need to show the statement when $r \geq 1$. Now, if $n < s$, then both sides are $0$, thanks to our convention that $\xi_{r-e}(m) = 0$ for any $m < 0$. Hence, we only need to show the statement when $r \geq 1$ and $n \geq s$. For $n = s$, consider $r = 1$ first. The left-hand side is precisely the number of $A \in \Mat_{n}(\bF_{q})$ such that $f_{A}(t) = P(t)$, and this is equal to
\[\frac{|\GL_{n}(\bF_{q})|}{q^{n} - 1}\]

by Lemma \ref{matchar}, but note that this quantity is equal to the right-hand side. If $n = s$ and $r \geq 2$, then the left-hand side is $0$. On the right-hand side, since $n = s$, each term has a factor of the form $\xi_{r-e}(s - es) = \xi_{r-e}(s(1-e))$ for $1 \leq e \leq s$, so every term except when $e = 1$ is $0$. Moreover, when $e = 1$, we have the factor $\xi_{r-1}(0)$, which is $0$ because $r \geq 2$ forces $r-1 \geq 1$, so we establish the identity whenever $n = s$. Thus, we only need to show the statement for $r \geq 1$ and $n > s$. 

\hspace{3mm} For $r \geq 1$ and $n > s = \deg(P)$, applying Lemma \ref{matchar}, we have
\begin{align*}
\xi_{r}(n) - \psi_{r}(n; P) &= \#\{A \in \Mat_{n}(\bF_{q}) : f_{A}(t) \text{ is square-free with $r$ irreducible factors such that } P(t) | f_{A}(t)\}\\
&= \sum_{\substack{P_{1}, \dots, P_{r-1} \in \mc{P}_{q} \sm \{P\} \text{ distinct:} \\ \deg(P_{1}) + \cdots + \deg(P_{r-1}) = n-s}} \#\{A \in \Mat_{n}(\bF_{q}) : f_{A}(t) = P_{1}(t) \cdots P_{r-1}(t)P(t)\} \\
&= \sum_{\substack{P_{1}, \dots, P_{r-1} \in \mc{P}_{q} \sm \{P\} \text{ distinct:} \\ \deg(P_{1}) + \cdots + \deg(P_{r-1}) = n-s}}\frac{|\GL_{n}(\bF_{q})|}{(q^{\deg(P_{1})} - 1) \cdots (q^{\deg(P_{r-1})} - 1)(q^{s} - 1)} \\
&= \frac{|\GL_{n}(\bF_{q})|}{(q^{s} - 1)|\GL_{n-s}(\bF_{q})|}\sum_{\substack{P_{1}, \dots, P_{r-1} \in |\bA^{n}(\bF_{q})| \sm \{P\} \text{ distinct:} \\ \deg(P_{1}) + \cdots + \deg(P_{r-1}) = n-s}}\frac{|\GL_{n-s}(\bF_{q})|}{(q^{\deg(P_{1})} - 1) \cdots (q^{\deg(P_{r-1})} - 1)} \\
&=  \frac{|\GL_{n}(\bF_{q})|}{(q^{s} - 1)|\GL_{n-s}(\bF_{q})|}\psi_{r-1}(n-s; P).
\end{align*}

So far, we have
\begin{equation}\label{aux}
\xi_{r}(n) = \psi_{r}(n; P) + \frac{|\GL_{n}(\bF_{q})|}{(q^{s} - 1)|\GL_{n-s}(\bF_{q})|}\psi_{r-1}(n-s; P).
\end{equation}

We showed the identity (\ref{aux}) for $r \geq 1$ and $n > s$. However, note that this identity still holds when $r \geq 1$ and $n = s$. When $n < s$, we have $\xi_{r}(n) = \psi_{r}(n; P)$ and since $\psi_{r-1}(n-s; P) = 0$, the identity (\ref{aux}) still holds with our convention $|\GL_{n-s}(\bF_{q})| = 1$. This tells us that the identity (\ref{aux}) holds for any $r \in \bZ_{\geq 1}$ and $n \in \bZ$. We have $\xi_{0}(0) = 1 = \psi_{0}(0; P)$ and  $\xi_{0}(n) = 0 = \psi_{0}(n; P)$ whenever $n \neq 0$. If $r < 0$, then $\xi_{r}(n) = 0 = \psi_{r}(n; P)$, so the identity (\ref{aux}) holds for any integers $r$ and $n$ without any restrictions.

\hspace{3mm} Going back to proving the desired statement, we are in the case $r \geq 1$ and $n > s$. We have the following $r$ identities by using the identity (\ref{aux}) with different inputs:
\bi
	\item $\xi_{r-1}(n-s) = \psi_{r-1}(n-s; P) + \frac{|\GL_{n-s}(\bF_{q})|}{(q^{s} - 1)|\GL_{n-2s}(\bF_{q})|}\psi_{r-2}(n-2s; P),$
	\item $\xi_{r-2}(n-2s) = \psi_{r-2}(n-2s; P) + \frac{|\GL_{n-2s}(\bF_{q})|}{(q^{s} - 1)|\GL_{n-3s}(\bF_{q})|}\psi_{r-3}(n-3s; P),$ \\
	\item $\xi_{r-3}(n-3s) = \psi_{r-3}(n-3s; P) + \frac{|\GL_{n-3s}(\bF_{q})|}{(q^{s} - 1)|\GL_{n-4s}(\bF_{q})|}\psi_{r-4}(n-4s; P),$ \\
	$\vdots$ \\
	\item $\xi_{1}(n-(r-1)s) = \psi_{1}(n-(r-1)s ; P) + \frac{|\GL_{n-(r-1)s}(\bF_{q})|}{(q^{s} - 1)|\GL_{n-rs}(\bF_{q})|}\psi_{0}(n-rs; P),$
	\item $\xi_{0}(n-rs) = \psi_{0}(n-rs ; P).$
\ei

Multiplying suitable quantities to the both sides of each identity, we get
\bi
	\item $\frac{|\GL_{n}(\bF_{q})|}{(q^{s} - 1)|\GL_{n-s}(\bF_{q})|}\xi_{r-1}(n-s) \\
	= \frac{|\GL_{n}(\bF_{q})|}{(q^{s} - 1)|\GL_{n-s}(\bF_{q})|}\psi_{r-1}(n-s; P) + \frac{|\GL_{n}(\bF_{q})|}{(q^{s} - 1)^{2}|\GL_{n-2s}(\bF_{q})|}\psi_{r-2}(n-2s; P),$
	\item $\frac{|\GL_{n}(\bF_{q})|}{(q^{s} - 1)^{2}|\GL_{n-2s}(\bF_{q})|}\xi_{r-2}(n-2s) \\
	= \frac{|\GL_{n}(\bF_{q})|}{(q^{s} - 1)^{2}|\GL_{n-2s}(\bF_{q})|}\psi_{r-2}(n-2s; P) + \frac{|\GL_{n}(\bF_{q})|}{(q^{s} - 1)^{3}|\GL_{n-3s}(\bF_{q})|}\psi_{r-3}(n-3s; P),$
	\item $ \frac{|\GL_{n}(\bF_{q})|}{(q^{s} - 1)^{3}|\GL_{n-3s}(\bF_{q})|}\xi_{r-3}(n-3s) \\
	= \frac{|\GL_{n}(\bF_{q})|}{(q^{s} - 1)^{3}|\GL_{n-3s}(\bF_{q})|}\psi_{r-3}(n-3s; P) + \frac{|\GL_{n}(\bF_{q})|}{(q^{s} - 1)^{4}|\GL_{n-4s}(\bF_{q})|}\psi_{r-4}(n-4s; P),$ \\
	$\vdots$ \\
	\item $\frac{|\GL_{n}(\bF_{q})|}{(q^{s} - 1)^{r-1} |\GL_{n-(r-1)s}(\bF_{q})|}\xi_{1}(n-(r-1)s) \\= \frac{|\GL_{n}(\bF_{q})|}{(q^{s} - 1)^{r-1} |\GL_{n-(r-1)s}(\bF_{q})|}\psi_{1}(n-(r-1)s ; P) + \frac{|\GL_{n}(\bF_{q})|}{(q^{s} - 1)^{r} |\GL_{n-rs}(\bF_{q})|}\psi_{0}(n-rs; P),$
	\item $\frac{|\GL_{n}(\bF_{q})|}{(q^{s} - 1)^{r}|\GL_{n-rs}(\bF_{q})|}\xi_{0}(n-rs) = \frac{|\GL_{n}(\bF_{q})|}{(q^{s} - 1)^{r} |\GL_{n-rs}(\bF_{q})|}\psi_{0}(n-rs ; P).$
\ei

Then we note that alternating the right-hand sides of the above identities yields cancellations to give us

\begin{align*}
\xi_{r}(n) - \psi_{r}(n; P) &= \frac{|\GL_{n}(\bF_{q})|}{(q^{s} - 1)|\GL_{n-s}(\bF_{q})|}\psi_{r-1}(n-s; P) \\
&= |\GL_{n}(\bF_{q})|\sum_{e=1}^{r}(-1)^{e-1}\frac{\xi_{r-e}(n-es)}{(q^{s}-1)^{e} |\GL_{n-es}(\bF_{q})|},
\end{align*}

as desired.
\end{proof}

\hspace{3mm} We now give our proof of Theorem \ref{recursion}:

\begin{proof}[Proof of Theorem \ref{recursion}] We assume $n, r \geq 1$ because otherwise the result is trivial. Given any $P \in \mc{P}_{q}$, we recall that 
\[\xi_{r}(n) - \psi_{r}(n ; P) = \#\left\{\begin{array}{c}
A \in \Mat_{n}(\bF_{q}) : 
f_{A}(t) \text{ is square-free with} \\
r\text{ irreducible factors and } P(t) | f_{A}(t)
\end{array}
\right\}.\]

This implies that
\[\prod_{\substack{A \in \Mat_{n}(\bF_{q}): \\ f_{A}(t) \text{ square-free with} \\ r \text{ irreducible factors}}}f_{A}(t) = \prod_{s=1}^{n}\prod_{\substack{P \in \mc{P}_{q}: \\ \deg(P) = s}}P(t)^{\xi_{r}(n) - \psi_{r}(n ; P)}.\]

Taking degrees of both sides, we get

\begin{align*}
n\xi_{r}(n) &= \sum_{s=1}^{n}s\pi(s)(\xi_{r}(n) - \psi_{r}(n;P)) \\
&= |\GL_{n}(\bF_{q})|\sum_{s=1}^{n}\sum_{e=1}^{r}(-1)^{e-1} \frac{s\pi(s)\xi_{r-e}(n-es)}{(q^{s}-1)^{e} |\GL_{n - es}(\bF_{q})|},
\end{align*}

where we used Lemma \ref{matchar} for the last equality. This finishes the proof. 
\end{proof}

\

\subsection{Non-square-free case} In this subsection, we provide a recursive formula for
\[\xi'_{r}(n) := \#\left\{\begin{array}{c}
A \in \Mat_{n}(\bF_{q}) : 
f_{A}(t) \text{ has } r\text{ irreducible factors} \\
\text{counting with multiplicity}
\end{array}
\right\}.\]

As in the square-free case, it turns out to be more convenient to find a recursive formula for
\[\phi'_{r}(n) := \xi'_{r}(n) / |\GL_{n}(\bF_{q})|.\]

Given any $P \in \mc{P}_{q}$, we define
\[\psi'_{r}(n; P) := \#\left\{\begin{array}{c}
A \in \Mat_{n}(\bF_{q}) : 
f_{A}(t) \text{ has } r\text{ irreducible factors} \\
\text{counting with multiplicity and } P(t) \nmid f_{A}(t)
\end{array}
\right\},\]

and for $e \in \bZ_{\geq 0}$, we define 
\[\psi'_{r}(n; P; e) := \#\left\{\begin{array}{c}
A \in \Mat_{n}(\bF_{q}) : 
f_{A}(t) \text{ has } r\text{ irreducible factors} \\
\text{counting with multiplicity while} \\
P(t)^{e} | f_{A}(t) \text{ and } P(t)^{e+1} \nmid f_{A}(t)
\end{array}
\right\}.\]

Note that $\psi'_{r}(n; P; 0) = \psi'_{r}(n; P)$. Define 
\[\alpha_{r}(n; P; e) := \psi'_{r}(n; P; e)/|\GL_{n}(\bF_{q})|\]

and $\alpha_{r}(n; P) := \alpha_{r}(n; P; 0)$.

\begin{rmk} We note that $\xi'_{0}(0) = \phi'_{0}(0) = 1$. If $r \geq 1$, then $\xi'_{r}(0) = \phi'_{r}(0) = 0$, and if $n \geq 1$, then $\xi'_{0}(n) = \phi'_{0}(n) = 0$. If $n < 0$ or $r < 0$, we set $\xi'_{r}(n) = \phi'_{r}(n) := 0$ and $\psi_{r}(n; P) := 0$. We also set $\psi'_{r}(n; P; e) := \alpha_{r}(n; P; e) = 0$ if $n < e\deg(P)$ or $r < e$.
\end{rmk}

\begin{lem}\label{alc} Fix $n, r, s \in \bZ_{\geq 1}$.  Let $P \in \mc{P}_{q}$ such that $\deg(P) = s$. For any $e \in \bZ_{\geq 0}$, we have
\[\alpha_{r}(n; P; e) = c_{e,s}\alpha_{r-e}(n - es; P),\]

where
\[c_{e,s} := \lt( q^{es}\prod_{i=1}^{e}(1 - q^{-is}) \rt)^{-1}.\]
\end{lem}

\begin{proof} If $n < r$, both sides are $0$, so assume $n \geq r$. If $r < e$, then both sides are $0$, so assume $r \geq e$. If $n < es$, then both sides are $0$, so assume $n \geq es$. Then we have
\begin{align*}
\psi'_{r}(n; P; e) &= \sum_{\substack{P_{1}, \dots, P_{l} \in \mc{P}_{q} \sm \{P\}, \\ e_{1}, \dots, e_{l} \in \bZ_{\geq 1}: \\ e_{1} + \cdots + e_{l} = r - e \text{ and} \\ e_{1}\deg(P_{1}) + \cdots + e_{l}\deg(P_{l}) = n - es}} \#\{A \in \Mat_{n}(\bF_{q}) : f_{A}(t) = P_{1}(t)^{e_{1}} \cdots P_{l}(t)^{e_{l}}P(t)^{e}\} \\
&= \sum_{\substack{P_{1}, \dots, P_{l} \in \mc{P}_{q} \sm \{P\}, \\ e_{1}, \dots, e_{l} \in \bZ_{\geq 1}: \\ e_{1} + \cdots + e_{l} = r - e \text{ and} \\ e_{1}\deg(P_{1}) + \cdots + e_{l}\deg(P_{l}) = n - es}} \frac{|\GL_{n}(\bF_{q})|}{q^{n} (1 - q^{-s})(1 - q^{-2s}) \cdots (1 - q^{-es}) \prod_{i=1}^{l} \prod_{j=1}^{e_{i}}(1 - q^{-j\deg(P_{i})})} \\
&= \frac{c_{e,s}|\GL_{n}(\bF_{q})|}{|\GL_{n-es}(\bF_{q})|}\sum_{\substack{P_{1}, \dots, P_{l} \in \mc{P}_{q} \sm \{P\}, \\ e_{1}, \dots, e_{l} \in \bZ_{\geq 1}: \\ e_{1} + \cdots + e_{l} = r - e \text{ and} \\ e_{1}\deg(P_{1}) + \cdots + e_{l}\deg(P_{l}) = n - es}} \frac{|\GL_{n-es}(\bF_{q})|}{q^{n-es} \prod_{i=1}^{l} \prod_{j=1}^{e_{i}}(1 - q^{-j\deg(P_{i})})} \\
&= \frac{c_{e,s}|\GL_{n}(\bF_{q})|}{|\GL_{n-es}(\bF_{q})|} \psi'_{r-e}(n-es; P),
\end{align*}

which gives the result.
\end{proof}

\begin{lem}\label{matid}  Fix $n, r, s \in \bZ_{\geq 1}$.  Let $P \in \mc{P}_{q}$ such that $\deg(P) = s$. Then
\[\phi'_{r-k}(n-ks) = \sum_{e=0}^{r-k}c_{e,s} \alpha_{r-(k+e)}(n - (k+e)s ; P)\]

for $0 \leq k \leq r$. In other words, we have
\[\bbm
\phi'_{r}(n) \\
\phi'_{r-1}(n-s) \\
\vdots \\
\phi'_{1}(n-(r-1)s) \\
\phi'_{0}(n-rs)
\ebm
=\bbm
c_{0,s} & c_{1,s} & \cdots & c_{r-1,s} & c_{r,s} \\
0 & c_{0,s} & \cdots & c_{r-2,s} & c_{r-1,s} \\
\vdots & \vdots & \cdots & \vdots & \vdots \\
0 & 0 & \cdots & c_{0,s} & c_{1,s} \\
0 & 0 & \cdots & 0 & c_{0,s}
\ebm
\bbm
\alpha_{r}(n; P) \\
\alpha_{r-1}(n-s; P) \\
\vdots \\
\alpha_{1}(n-(r-1)s; P) \\
\alpha_{0}(n-rs; P)
\ebm.\]
\end{lem}

\begin{proof} Fix $0 \leq k \leq r$. If $r < k$ or $n < ks$, then both sides are $0$, so assume that $r \geq k$ and $n \geq ks$. For any $A \in \Mat_{n-ks}(\bF_{q})$, there is a unique $0 \leq e \leq r-k$ such that $P(t)^{e} | f_{A}(t)$ and $P(t)^{e+1} \nmid f_{A}(t)$. This implies that
\[\xi'_{r-k}(n-ks) = \sum_{e=0}^{r-k}\psi'_{r-k}(n-ks; P; e),\]

so dividing by $|\GL_{n-ks}(\bF_{q})|$, we have
\[\phi'_{r-k}(n-ks) = \sum_{e=0}^{r-k}\al_{r-k}(n-ks; P; e).\]
\end{proof}

\begin{rmk} From now on, we write
\[E_{r,s} := \bbm
c_{0,s} & c_{1,s} & \cdots & c_{r-1,s} & c_{r,s} \\
0 & c_{0,s} & \cdots & c_{r-2,s} & c_{r-1,s} \\
\vdots & \vdots & \cdots & \vdots & \vdots \\
0 & 0 & \cdots & c_{0,s} & c_{1,s} \\
0 & 0 & \cdots & 0 & c_{0,s}
\ebm,\]

which still depends on $q$, but we drop the notation because we fix $q$. Recall that
\[c_{e,s} = \lt( q^{es}\prod_{i=1}^{e}(1 - q^{-is}) \rt)^{-1},\]

so in particular, we have $c_{0,s} = 1$. Hence, the matrix $E_{r,s}$ is invertible over $\bQ$ (i.e., is an element of $\GL_{r+1}(\bQ)$) and its inverse $E_{r,s}^{-1}$ is an upper-triangular matrix whose diagonal entries are $1$. 

\hspace{3mm} Given any matrix $A$ (of any size over any commutative ring), we denote by $A_{i,j}$ the $(i,j)$-entry of $A$. We note that applying $E_{r,s}^{-1}$ both sides on Lemma \ref{matid}, for $0 \leq e \leq r$, we get
\begin{equation}\label{al}
\al_{r-e}(n-es; P) = \sum_{j=e}^{r}(E_{r,s}^{-1})_{e+1,j+1} \phi'_{r-j}(n - js),
\end{equation}

which depends on $s = \deg(P)$ rather than $P$ itself.
\end{rmk}

\hspace{3mm} We now give a recursion formula for $\phi'_{r}(n)$:

\begin{thm}\label{recursion2} Given $n, r \in \bZ_{\geq 0}$, we have
\[n\phi'_{r}(n) = \sum_{s=1}^{n}\sum_{e=1}^{r}\sum_{j=e}^{r} s\pi(s)e c_{e,s}(E_{r,s}^{-1})_{e+1,j+1}\phi'_{r-j}(n-js).\]
\end{thm}

\begin{proof} We have
\[\prod_{\substack{A \in \Mat_{n}(\bF_{q}): \\ f_{A}(t) \text{ has } r \text{ irreducible factors} \\ \text{counting with multiplicity}}}f_{A}(t) = \prod_{s=1}^{n}\prod_{e=1}^{r}\prod_{\substack{P \in \mc{P}_{q}: \\ \deg(P) = s}}P(t)^{e\psi'_{r}(n ; P; e)},\]

so taking the degrees for both sides, we get
\[n\xi'_{r}(n) = \sum_{s=1}^{n}\sum_{e=1}^{r}s\pi(s) e \psi'_{r}(n; P; e),\]

so dividing by $|\GL_{n}(\bF_{q})|$ and applying Lemma \ref{alc}, we have
\begin{align*}
n\phi'_{r}(n) &= \sum_{s=1}^{n}\sum_{e=1}^{r}s\pi(s) e \alpha_{r}(n; P; e) \\
&= \sum_{s=1}^{n}\sum_{e=1}^{r}s\pi(s) e c_{e,s}\alpha_{r-e}(n - es; P) \\
&= \sum_{s=1}^{n}\sum_{e=1}^{r}\sum_{j=e}^{r}s\pi(s) e c_{e,s}(E_{r,s}^{-1})_{e+1,j+1} \phi'_{r-j}(n - js),
\end{align*}

where we used (\ref{al}) for the last equality.
\end{proof}

\


\section{Proof of Theorem \ref{main}}

\subsection{Square-free case} Using the notation in Theorem \ref{recursion}, we note that the square-free case of Theorem \ref{main} is equivalent to the following:

\begin{thm}\label{main2} Let $r \in \bZ_{\geq 1}$. Then
\[\phi_{r}(n) = \frac{(\log n)^{r-1}}{(r-1)!n} + O\left( \frac{(\log n)^{r-2}}{n} \right),\]

\

where the implied constant may depend on $r$ but not on $q$ nor $n$.
\end{thm}

\begin{rmk} Theorem \ref{main2} immediately implies that 
\begin{align*}
&\underset{A \in \Mat_{n}(\bF_{q})}{\Prob}
\left(\begin{array}{c}
f_{A}(t) \text{ is square-free with} \\
r\text{ irreducible factors}
\end{array}\right) = \left(\frac{(\log n)^{r-1}}{(r-1)!n} + O\left( \frac{(\log n)^{r-2}}{n} \right)\right)\prod_{i=1}^{n}(1 - q^{-i}).
\end{align*}

Together with Lemma \ref{rec}, Theorem \ref{main2} also implies that
\[\#\left\{\begin{array}{c}
A \in \Mat_{n}(\bF_{q}) : 
f_{A}(t) \text{ is square-free with} \\
r\text{ irreducible factors and } A \notin \GL_{n}(\bF_{q})
\end{array}
\right\} = |\GL_{n}(\bF_{q})| \cdot O\left( \frac{(\log n)^{r-2}}{n} \right)\]

where the implied constant may depend on $r$ but not on $q$ nor $n$. (Note that $t | f_{A}(t)$ if and only if $A \notin \GL_{n}(\bF_{q})$.) As a result, we have
\begin{align*}
&\underset{A \in \GL_{n}(\bF_{q})}{\Prob}
\left(\begin{array}{c}
f_{A}(t) \text{ is square-free with} \\
r\text{ irreducible factors}
\end{array}\right) = \frac{(\log n)^{r-1}}{(r-1)!n} + O\left( \frac{(\log n)^{r-2}}{n} \right).
\end{align*}
\end{rmk}

\subsection{Non-square-free case} We note that the non-square-free case of Theorem \ref{main} is equivalent to the following:

\begin{thm}\label{main3} Let $r \in \bZ_{\geq 1}$. Then
\[\phi'_{r}(n) = \frac{(\log n)^{r-1}}{(r-1)!n} + O\left( \frac{(\log n)^{r-2}}{n} \right),\]

\

where the implied constant may depend on $r$ but not on $q$ nor $n$.
\end{thm}

\begin{rmk} Theorem \ref{main3} immediately implies that 
\begin{align*}
&\underset{A \in \Mat_{n}(\bF_{q})}{\Prob}
\left(\begin{array}{c}
f_{A}(t) \text{ has } r \text{ irreducible factors} \\
\text{counting with multiplicity}
\end{array}\right) = \left(\frac{(\log n)^{r-1}}{(r-1)!n} + O\left( \frac{(\log n)^{r-2}}{n} \right)\right)\prod_{i=1}^{n}(1 - q^{-i}).
\end{align*}

Applying Lemma \ref{alc}, Theorem \ref{main3}, and the fact that $|(E_{r,s}^{-1})_{i,j}| = O(1)$ for any $1 \leq i, j \leq r+1$ whose implied constant only depends on $r$ (which follows from Lemma \ref{ij}), we get
\[\#\left\{\begin{array}{c}
A \in \Mat_{n}(\bF_{q}) : 
f_{A}(t) \text{ has } r \text{ irreducible factors} \\
\text{counting with multiplicity and } A \notin \GL_{n}(\bF_{q})
\end{array}
\right\} = |\GL_{n}(\bF_{q})| \cdot O\left( \frac{(\log n)^{r-2}}{n} \right).\]

As a result, we have
\begin{align*}
&\underset{A \in \GL_{n}(\bF_{q})}{\Prob}
\left(\begin{array}{c}
f_{A}(t) \text{ has } r \text{ irreducible factors}\\
\text{counting with multiplicity}
\end{array}\right) = \frac{(\log n)^{r-1}}{(r-1)!n} + O\left( \frac{(\log n)^{r-2}}{n} \right).
\end{align*}
\end{rmk}

\

\subsection{Useful lemmas} We provide useful lemmas for our proof of Theorems \ref{main2} and \ref{main3} in this subsection. 

\begin{notn} Given $d_{1}, \dots, d_{r} \in \bZ_{\geq 1}$, denote by $M_{q}(d_{1}, \dots, d_{r})$ be the number of monic polynomials in $\bF_{q}[t]$ of the form $P_{1}(t) \cdots P_{r}(t)$, where $P_{1}(t), \dots, P_{r}(t)$ are monic irreducible polynomials with $\deg(P_{i}) = d_{i}$ for $1 \leq i \leq r$.
\end{notn}

\hspace{3mm} The following lemma is a reformulation of Cohen's result mentioned in the introduction:

\begin{lem}[Cohen]\label{Cohen} Let $n, r \in \bZ_{\geq 1}$. We have
\[\frac{1}{q^{n}}\sum_{\substack{d_{1} \geq \cdots \geq d_{r} \geq 1 \\ d_{1} + \cdots + d_{r} = n}}M_{q}(d_{1}, \dots, d_{r}) = \frac{(\log n)^{r-1}}{(r-1)! n} + O\left(\frac{(\log n)^{r-2}}{n}\right),\]

where the implied constant may depend on $r$ but does not depend on $q$ nor $n$.
\end{lem}

\hspace{3mm} The following lemma follows from \cite[(3.12) and (3.13)]{Coh}, but the proof we include seems slightly different from Cohen's, so we include it here since this lemma is fundamental for our arguments later:

\begin{lem}\label{calc} Let $n, r \in \bZ_{\geq 1}$. Then
\[\sum_{s=1}^{n-1}\frac{(\log(s))^{r-1}}{s} = \frac{(\log(n))^{r}}{r} + O(1),\]

where the implied constant may depend on $r$ but not on $n$.
\end{lem}

\begin{proof} If $n = 1$, the sum is equal to $0$, so we may assume that $n \geq 2$. We note that
\[\int_{1}^{n} \frac{(\log(x))^{r-1}}{x} dx = \int_{0}^{\log(n)} y^{r-1} dy = \frac{(\log(n))^{r}}{r},\]

\

so we may use the integral on the left-hand side instead of the quantity on the right-hand side. Consider the real valued function $f(x) = (\log(x))^{r-1}/x$ for $x \geq 1$. Note that for $x \geq 1$, we have 
\bi
	\item $f'(x) \leq 0$ (i.e., $f(x)$ is decreasing) if and only if $x \geq e^{r-1}$;
	\item $f'(x) \geq 0$ (i.e., $f(x)$ is increasing) if and only if $x \leq e^{r-1}$.
\ei
	
Thus, if $n \leq e^{r-1}$, then we have
\[0 \leq \int_{1}^{n}f(x)dx - \sum_{s=1}^{n-1}f(s) \leq f(n) - f(1)\]

so that
\[0 \leq  \int_{1}^{n}\frac{(\log(x))^{r-1}}{x}dx - \sum_{s=1}^{n-1}\frac{(\log(s))^{r-1}}{s} \leq \frac{(\log(n))^{r-1}}{n} \leq \frac{(\log(e^{r-1}))^{r-1}}{e^{r-1}} = \frac{(r-1)^{r-1}}{e^{r-1}}.\]

If $n > e^{r-1}$, then
\begin{align*}
0 \leq \lt|\int_{1}^{n}\frac{(\log(x))^{r-1}}{x}dx - \sum_{s=1}^{n-1}\frac{(\log(s))^{r-1}}{s}\rt| &\leq 2\int_{c}^{c+1}\frac{(\log(x))^{r-1}}{x}dx \\
&= \frac{2}{r}\left((\log(c+1))^{r} - (\log(c))^{r}\right)
\end{align*}

where $c = \lfloor e^{r-1} \rfloor$. Since both bounds only depend on $r$ but not on $n$, this finishes the proof.
\end{proof}

\begin{lem}\label{ij} For any $r, s \in \bZ_{\geq 1}$ and $1 \leq i < j \leq r+1$, we have
\[|(E_{r,s}^{-1})_{i,j}| \leq O(q^{-(j-i)s}),\]

where the implied constant only depends on $r$.
\end{lem}

\begin{proof}  Let $N_{r,s} := E_{r,s} - I_{r+1}$. Since $N_{r,s}$ is an $(r+1) \times (r+1)$ matrix that is strictly upper-triangular, we have $N_{r,s}^{r+1} = 0$. Since $E_{r,s} = I_{r+1} + N_{r,s}$, we have
\[E_{r,s}^{-1} = \sum_{k=0}^{r}(-1)^{k}N_{r,s}^{k}.\]
Thus, it follows that
\[|(E_{r,s}^{-1})_{i,j}| \leq \sum_{k=0}^{r}(N_{r,s}^{k})_{i,j}.\]

Since $N_{r,s}^{k}$ is strictly upper-triangular, we have
\[(N_{r,s}^{k})_{i,j} = \sum_{i=i_{0} < i_{1} < \cdots < i_{k} = j}\prod_{u=1}^{k}(N_{r,s})_{i_{u-1},i_{u}}.\]

Now, note that $c_{e,s} = (q^{es}\prod_{i=1}^{e}(1-q^{-is}))^{-1} \leq (q^{es}\prod_{i=1}^{\infty}(1-2^{-i}))^{-1}  \leq 4q^{-es}$, so $(N_{r,s})_{i_{u-1},i_{u}} = c_{i_{u}-i_{u-1},s} \leq 4 q^{-(i_{u} - i_{u-1})s}$. This implies that
\[(N_{r,s}^{k})_{i,j} \leq \sum_{i=i_{0} < i_{1} < \cdots < i_{k} = j}\prod_{u=1}^{k} 4 q^{-(i_{u} - i_{u-1})s} = 4^{k}{j-i + 1 \choose k+1}q^{-(j-i)s}\]

for any $0 \leq k \leq r$, so the result follows.
\end{proof}

\

\subsection{Proof of Theorem \ref{main2}} In this subsection, we prove Theorem \ref{main2}. We start with the following lemma:

\begin{lem}\label{important} Let $n, r \in \bZ_{\geq 1}$ and $e \in \bZ_{\geq 2}$. Then
\[\sum_{s=1}^{n}\frac{s \pi(s)\phi_{r-e}(n-es)}{(q^{s} - 1)^{e}} = O(\log(n)^{r-3}),\]

where the implied constant may depend on $r$ and $e$ but not on $q$ nor $n$. In particular, by Theorem \ref{recursion}, we have
\[n\phi_{r}(n) = \lt(\sum_{s=1}^{n}\frac{s\pi(s)\phi_{r-1}(n-s)}{q^{s} - 1}\rt) + O(\log(n)^{r-3}).\]

where the implied constant may depend on $r$ but not on $q$ nor $n$.
\end{lem}

\begin{proof} If $r < e$, then the sum in the statement is $0$, so we may assume that $r \geq e$. Consider the case that $r = e \geq 2$. Then
\[\sum_{s=1}^{n}\frac{s \pi(s)\phi_{r-e}(n-es)}{(q^{s} - 1)^{e}} = \sum_{s=1}^{n}\frac{s \pi(s)\phi_{0}(n-rs)}{(q^{s} - 1)^{r}} = \left\{
	\begin{array}{ll}
	0 \mbox{ if } r \nmid n,  \\
	\frac{n\pi(n/r)}{r(q^{n/r}-1)^{r}} = O\left( \frac{1}{q^{n(r-1)/r}} \right)  \mbox{ if } r | n,
	\end{array}\right.\]
since we consider $r$ as a constant. The implied constant only depends on $r$ but not on $q$ nor $n$. Thus, the result immediately follows when $r = e$.

\hspace{3mm} It remains to consider the case $r > e$. By Lemma \ref{matchar}, we have
\[\xi_{r}(n) = \sum_{\substack{d_{1} \geq \cdots \geq d_{r} \geq 1 \\ d_{1} + \cdots + d_{r} = n}}\frac{M_{q}(d_{1}, \dots, d_{r})|\GL_{n}(\bF_{q})|}{(q^{d_{1}} -1) \cdots (q^{d_{r}} -1)},\]

where the notation $\xi_{r}(n)$ is as in Theorem \ref{recursion}. Thus, we have
\[\phi_{r}(n) = \frac{\xi_{r}(n)}{|\GL_{n}(\bF_{q})|}= \sum_{\substack{d_{1} \geq \cdots \geq d_{r} \geq 1 \\ d_{1} + \cdots + d_{r} = n}}\frac{M_{q}(d_{1}, \dots, d_{r})}{(q^{d_{1}} -1) \cdots (q^{d_{r}} -1)},\]

so noting that $\frac{1}{q^{d_{i}}-1} \leq \frac{2}{q^{d_{i}}}$, we have
\[ \frac{1}{q^{n}}\sum_{\substack{d_{1} \geq \cdots \geq d_{r} \geq 1 \\ d_{1} + \cdots + d_{r} = n}}M_{q}(d_{1}, \dots, d_{r}) \leq \phi_{r}(n) \leq  \frac{2^{r}}{q^{n}}\sum_{\substack{d_{1} \geq \cdots \geq d_{r} \geq 1 \\ d_{1} + \cdots + d_{r} = n}} M_{q}(d_{1}, \dots, d_{r}).\]

Applying Lemma \ref{Cohen}, we note that there are some positive real constants $c_{1} < c_{2}$ that only depend on $r$ (but not on $q$ nor $n$) such that
\[c_{1}\frac{(\log n)^{r-1}}{(r-1)! n} \leq \phi_{r}(n) \leq c_{2}\frac{(\log n)^{r-1}}{(r-1)! n}.\]

Thus, when $n - es \geq 1$, we have
\begin{align*}
\frac{s \pi(s)\phi_{r-e}(n-es)}{(q^{s} - 1)^{e}} &\leq \frac{c_{2}q^{s}(\log(n-es))^{r-e-1}}{(r-e-1)!(q^{s} - 1)^{e}(n-es)}\\
&\leq \frac{2^{e}c_{2}q^{s}(\log(n-es))^{r-e-1}}{(r-e-1)!q^{se}(n-es)}\\
&= \frac{2^{e}c_{2}(\log(n-es))^{r-e-1}}{(r-e-1)!q^{s(e-1)}(n-es)}.
\end{align*}

If $n - es \geq 1$, then $n \geq e$. Since $e \geq 2$, the above inequality gives
\[\frac{s \pi(s)\phi_{r-e}(n-es)}{(q^{s} - 1)^{e}} \leq \frac{2^{e}c_{2}(\log(n))^{r-3}}{q^{s}}.\]

Moreover, we note that the inequality above is true even without the condition $n - es \geq 1$ as the quantity on the left-hand side is $0$ when $n - es < 0$. We conclude that
\[\sum_{s=1}^{n}\frac{s \pi(s)\phi_{r-e}(n-es)}{(q^{s} - 1)^{e}} \leq 2^{e}c_{2}(\log(n))^{r-3}\frac{q^{-1}(1 - q^{-n})}{1 - q^{-1}} \leq 2^{e}c_{2}(\log(n))^{r-3},\]

as desired.
\end{proof}

\hspace{3mm} We are now ready to prove Theorem \ref{main2}:

\begin{proof}[Proof of Theorem \ref{main2}] We proceed by induction on $r$, the number of irreducible factors of $f_{A}(t)$ for a random matrix $A \in \Mat_{n}(\bF_{q})$. We have
\[\phi_{1}(n) = \frac{\pi(n)}{q^{n}-1}  = \frac{1}{q^{n} - 1}\left( \frac{q^{n}}{n} + O\lt(\frac{q^{n/2}}{n}\rt) \right) = \frac{1}{n} + O\left( \frac{1}{nq^{n/2}} \right),\]

where the implied constant does not depend on $q$ nor $n$. We note that this gives a stronger bound than the claimed one $O((n\log(n))^{-1})$, so this is sufficient for the case $r = 1$.

\hspace{3mm} Next, consider the case $r = 2$. By Theorem \ref{recursion}, we have
\begin{align*}
n\phi_{2}(n) &=  \sum_{s=1}^{n}\sum_{e=1}^{2}(-1)^{e-1} \frac{s\pi(s)\phi_{2-e}(n-es)}{(q^{s}-1)^{e}} \\
&= \sum_{s=1}^{n} \lt( \frac{s\pi(s)\phi_{1}(n-s)}{q^{s}-1} - \frac{s\pi(s)\phi_{0}(n-2s)}{(q^{s}-1)^{2}} \rt) \\
&= \sum_{s=1}^{n-1} \lt( \frac{s\pi(s)\phi_{1}(n-s)}{q^{s}-1} - \frac{s\pi(s)\phi_{0}(n-2s)}{(q^{s}-1)^{2}} \rt).
\end{align*}

We have 
\begin{equation}\label{est}
\frac{s\pi(s)}{q^{s}-1} = 1 + O\lt(\frac{1}{q^{s/2}}\rt),
\end{equation}

where the implied constant does not depend on $q$ or $n$. Since
\[\phi_{1}(n-s) = \frac{1}{n-s} + O\lt( \frac{1}{(n-s)q^{(n-s)/2}} \rt),\]

the above recursive formula with Lemma \ref{calc} with $r = 1$ implies that
\[n\phi_{2}(n) = \log n + O(1),\]

where the implied constant does not depend on $q$ nor $n$. This proves the theorem for the case $r = 2$.

\hspace{3mm} For the induction hypothesis, we assume that the statement holds for $\phi_{1}(n), \phi_{2}(n), \dots, \phi_{r-1}(n)$ with $r-1 \geq 2$. We then show the statement for $\phi_{r}(n)$. By the second statement of Lemma \ref{important}, we have
\begin{align*}
n\phi_{r}(n) &= \lt(\sum_{s=1}^{n}\frac{s\pi(s)\phi_{r-1}(n-s)}{q^{s} - 1}\rt) + O\lt((\log n)^{r-3}\rt)\\
&= \lt(\sum_{s=1}^{n-1}\frac{s\pi(s)\phi_{r-1}(n-s)}{q^{s} - 1}\rt) + O\lt((\log n)^{r-3}\rt),
\end{align*}

where the second identity used the fact that we are focusing on $r \geq 2$ (in fact, we are assuming $r \geq 3$) so that $\phi_{r-1}(0) = 0$. On top of this, we use (\ref{est}) and induction hypothesis:
\[\phi_{r-1}(n-s) = \frac{(\log(n - s))^{r-2}}{(r-2)!(n-s)} + O\lt( \frac{(\log(n-s))^{r-3}}{n-s} \rt).\]

Thus, we have
\begin{align*}
n\phi_{r}(n) &= O\lt(\log(n)^{r-3}\rt) + \sum_{s=1}^{n-1} \lt(1 + O\lt(\frac{1}{q^{s/2}}\rt)\rt) \lt(\frac{(\log(n - s))^{r-2}}{(r-2)!(n-s)} + O\lt( \frac{(\log(n-s))^{r-3}}{n-s} \rt)\rt) \\
&= O\lt(\log(n)^{r-3}\rt) + \lt( \sum_{s=1}^{n-1} \frac{(\log(n - s))^{r-2}}{(r-2)!(n-s)} \rt) +  O\lt(\sum_{s=1}^{n-1} \frac{(\log(n - s))^{r-2}}{(r-2)!q^{s/2}(n-s)}  \rt) + O\lt( \sum_{s=1}^{n-1} \frac{(\log(n-s))^{r-3}}{n-s} \rt) \\
&= O\lt(\log(n)^{r-3}\rt) + \frac{(\log(n))^{r-1}}{(r-1)!} +  O\lt( \sum_{s=1}^{n-1} \frac{(\log(n - s))^{r-2}}{(r-2)!q^{s/2}(n-s)}  \rt) + O\lt( (\log(n))^{r-2} \rt),
\end{align*}

where the last estimate uses Lemma \ref{calc}. Then note that the function $x \mapsto (\log(x))^{r-2}/x$ defined for $x \in [1, \infty)$ takes its global maximum at $x = e^{r-2}$, so
\begin{align*}
\sum_{s=1}^{n-1} \frac{(\log(n - s))^{r-2}}{(r-2)!q^{s/2}(n-s)} &= \frac{1}{(r-2)!}\sum_{s=1}^{n-1} \frac{(\log(n - s))^{r-2}}{q^{s/2}(n-s)} \\
&\leq \frac{1}{(r-2)!}\sum_{s=1}^{n-1} \frac{(\log(e^{r-2}))^{r-2}}{q^{s/2}e^{r-2}} \\
&= \frac{(r-2)^{r-2}}{(r-2)!e^{r-2}}\sum_{s=1}^{n-1} \frac{1}{q^{s/2}} \\
&= \frac{(r-2)^{r-2}}{(r-2)!e^{r-2}}\frac{q^{-1/2}(1 - q^{-(n-1)/2})}{1 - q^{-1/2}} \\
&\leq \frac{(r-2)^{r-2}}{(r-2)!e^{r-2}}\frac{q^{-1/2}}{1 - q^{-1/2}} \\
&\leq \frac{2(r-2)^{r-2}}{(r-2)!e^{r-2}}.
\end{align*}

This finishes the proof.
\end{proof}

\

\subsection{Proof of Theorem \ref{main3}} In this subsection, we prove Theorem \ref{main3}.

\begin{proof}[Proof of Theorem \ref{main3}] We proceed by induction on $r \in \bZ_{\geq 1}$. When $r = 1$, we have\footnote{All the big-O terms in this proof only depend on $r$.}
\[\phi'_{1}(n) = \phi_{1}(n) = \frac{1}{n}+ O\lt ( \frac{(\log n)^{-1}}{n} \rt)\]

by Theorem \ref{main2}. When $r=2$, Theorem \ref{recursion2} implies that
\begin{align*}
n\phi'_{2}(n) &= \sum_{s=1}^{n}\sum_{e=1}^{2}\sum_{j=e}^{2} s\pi(s)e c_{e,s}(E_{2,s}^{-1})_{e+1,j+1}\phi'_{2-j}(n-js) \\
&= \sum_{s=1}^{n}\lt(2s\pi(s) c_{2,s}(E_{2,s}^{-1})_{3,3}\phi'_{0}(n-2s) + \sum_{j=1}^{2} s\pi(s) c_{1,s}(E_{2,s}^{-1})_{2,j+1}\phi'_{2-j}(n-js) \rt) \\
&= \sum_{s=1}^{n} s\pi(s) (2c_{2,s} \phi'_{0}(n-2s) +  c_{1,s} \phi'_{1}(n-s) + c_{1,s}(E_{2,s}^{-1})_{2,3}\phi'_{0}(n-2s)) \\
&= \sum_{s=1}^{n-1} s\pi(s) (2c_{2,s} \phi'_{0}(n-2s) +  c_{1,s} \phi'_{1}(n-s) + c_{1,s}(E_{2,s}^{-1})_{2,3}\phi'_{0}(n-2s)).
\end{align*}

We know $\phi'_{0}(n-2s) = 0$ unless $s = n/2$ which sets $\phi'_{0}(n-2s) = 1$. We recall that
\[\frac{s\pi(s)}{q^{s}-1} = 1 + O(q^{-s/2}),\]

and we also must recall that
\[\phi'_{1}(n-s) = \phi_{1}(n-s) = \frac{1}{n-s} + O\lt( \frac{1}{(n-s)q^{(n-s)/2}} \rt)\]

in the proof of Theorem \ref{main2}. (Note that if we use $O(((n-s)\log(n-s))^{-1})$ instead, our argument below does not work.) Bounding the term $(E_{2,s}^{-1})_{2,3}$ using Lemma \ref{ij}, we have
\begin{align*}
n \phi'_{2}(n) &= O(1) + \sum_{s=1}^{n-1}s\pi(s)c_{1,s}\phi'_{1}(n-s) \\
&= O(1) + \sum_{s=1}^{n-1}\frac{s\pi(s)}{q^{s}-1}\lt( \frac{1}{n-s} + O\lt( \frac{1}{(n-s)q^{(n-s)/2}} \rt) \rt) \\
&= O(1) + \sum_{s=1}^{n-1}(1 + O(q^{-s/2}))\lt( \frac{1}{n-s} + O\lt( \frac{1}{(n-s)q^{(n-s)/2}} \rt) \rt) \\
&= O(1) + \lt(\sum_{s=1}^{n-1} \frac{1}{n-s}\rt) + O\lt(\sum_{s=1}^{n-1}\frac{q^{-s/2}}{n-s}\rt) + O\lt(\sum_{s=1}^{n-1}\frac{1}{(n-s)q^{(n-s)/2}}\rt) + O\lt(\sum_{s=1}^{n-1}\frac{1}{(n-s)q^{n/2}}\rt) \\
&= \log n + O(1)
\end{align*}

because
\[\sum_{s=1}^{n-1}q^{-s/2} \leq \sum_{s=1}^{\infty}2^{-s/2} = \frac{2^{-1/2}}{1 - 2^{-1/2}}.\]

Hence, we have
\[\phi'_{2}(n) = \frac{\log n}{n} + O(1/n),\]

proving Theorem \ref{main3} for the case $r = 2$.

\hspace{3mm} For the induction hypothesis, we assume that the result holds for $\phi'_{1}(n), \phi'_{2}(n), \dots, \phi'_{r-1}(n)$, and wish to show the result for $\phi'_{r}(n)$, where $r \geq 3$. By Theorem \ref{recursion2}, we have
\[n \phi'_{1}(n) = A_{n,r} + \sum_{j=2}^{r}\sum_{e=1}^{j}eB_{n,r,j,e}\]

where
\[A_{n,r} := \sum_{s=1}^{n}s\pi(s) c_{1,s}(E_{r,s}^{-1})_{2,2} \phi'_{r-1}(n-s) = \sum_{s=1}^{n}\frac{s\pi(s) \phi'_{r-1}(n-s)}{q^{s}-1}\]

and
\[B_{n,r,j,e} := \sum_{s=1}^{n} s\pi(s) c_{e,s}(E_{r,s}^{-1})_{e+1,j+1} \phi'_{r-j}(n-js).\]

We immediately note that for both $A_{n,r}$ and $B_{n,r,j,e}$ putting $s = n$ yields a zero summand. Thus, we have
\[A_{n,r} = \sum_{s=1}^{n-1}\frac{s\pi(s) \phi'_{r-1}(n-s)}{q^{s}-1}\]

and
\[B_{n,r,j,e} = \sum_{s=1}^{n-1} s\pi(s) c_{e,s}(E_{r,s}^{-1})_{e+1,j+1} \phi'_{r-j}(n-js).\]

We have
\[\frac{s\pi(s)}{q^{s}-1} = 1 + O(q^{-s/2}),\]

and the induction hypothesis implies that
\[\phi'_{r-1}(n-s) = \frac{(\log (n-s))^{r-2}}{(r-2)!(n-s)} + O\lt( \frac{(\log (n-s))^{r-3}}{n-s} \rt).\]

Hence, we have
\begin{align*}
A_{n,r} &= \sum_{s=1}^{n-1}\frac{s\pi(s) \phi'_{r-1}(n-s)}{q^{s}-1} \\
&= \sum_{s=1}^{n-1}(1 + O(q^{-s/2})) \lt(  \frac{(\log (n-s))^{r-2}}{(r-2)!(n-s)} + O\lt( \frac{(\log (n-s))^{r-3}}{n-s} \rt) \rt) \\
&= \lt(\sum_{s=1}^{n-1}\frac{(\log (n-s))^{r-2}}{(r-2)!(n-s)}\rt) + O\lt( \sum_{s=1}^{n-1} \frac{(\log (n-s))^{r-3}}{n-s} \rt) \\
&= \frac{(\log n)^{r-1}}{(r-1)!} + O\lt((\log n)^{r-2}\rt),
\end{align*}

where we used Lemma \ref{calc} and the inequality
\[\sum_{s=1}^{n-1} \frac{(\log(n - s))^{r-2}}{(r-2)!q^{s/2}(n-s)} \leq \frac{2(r-2)^{r-2}}{(r-2)!e^{r-2}},\]

which we established in our proof of Theorem \ref{main2}. From now on, we show
\[B_{n,r,j,e} = O((\log n)^{r-3})\]

for $2 \leq j \leq r$ and $1 \leq e \leq j$, where the implied constant only depends on $r$, which is enough to finish our proof.

\hspace{3mm} First, assume that $j = r$. Then
\[B_{n,r,j,e} = B_{n,r,r,e} = \sum_{s=1}^{n}s\pi(s)c_{e,s}(E_{r,s}^{-1})_{e+1,r+1} \phi'_{0}(n-rs) =\left\{
	\begin{array}{ll}
	0 \mbox{ if } r \nmid n,  \\
	(n/r)\pi(n/r) c_{e,n/r} (E_{r,n/r}^{-1})_{e+1,r+1} \mbox{ if } r | n.
	\end{array}\right.\]
	
If $r | n$, we have
\[(n/r)\pi(n/r) c_{e,n/r} = \frac{(n/r)\pi(n/r)}{q^{en/r}\prod_{i=1}^{e}(1 - q^{-in/r})} \leq \frac{(n/r)\pi(n/r)}{q^{n/r} - 1}  = O(1).\]

By Lemma \ref{ij}, we have
\[|(E_{r,n/r}^{-1})_{e+1,r+1}| = O(q^{-(r-e)n/r})\]

Thus, we have $B_{n,r,j,e} = O(1)$ when $j = r$. (Note that it is important that $r \geq 3$ here, which ensures that $1 \leq (\log n)^{r-3}$.)

\hspace{3mm} It remains to consider how to bound $B_{n,r,j,e}$ when $j < r$. (We also have $j \geq \max(e, 2)$.) Note that
\[\prod_{i=1}^{e}\frac{1}{1 - q^{-is}} \leq \prod_{i=1}^{\infty}\frac{1}{1 - 2^{-i}} \leq 4,\]

which implies that
\[c_{e,s} = \frac{1}{q^{es}\prod_{i=1}^{e}(1 - q^{-is})} = O(q^{-es}).\]

Since $s\pi(s) = O(q^{s})$ and $|(E_{r,s}^{-1})_{i,j}| = O(q^{-(j-i)s})$ for any $1 \leq i < j \leq r + 1$, we have
\[s\pi(s) c_{e,s}(E_{r,s}^{-1})_{e+1,j+1} = O(q^{s (1-e) })\]

Thus, we have
\begin{align*}
B_{n,r,j,e} &= \sum_{s=1}^{n-1} s\pi(s) c_{e,s}(E_{r,s}^{-1})_{e+1,j+1} \phi'_{r-j}(n-js) \\
&= O \lt( \sum_{1 \leq s < n/j} q^{s(1 - e)} \phi'_{r-j}(n-js) \rt) \\
&= O \lt( \sum_{1 \leq s < n/j} q^{s(1-e)}  \frac{(\log (n-js))^{r-j-1}}{n-js}\rt).
\end{align*}

This implies that when $e \geq 2$, we have
\[B_{n,r,j,e} =  O \lt( \sum_{1 \leq s < n/j} q^{-s}  \frac{(\log (n-js))^{r-j-1}}{n-js}\rt).\]

If $e = 1$, then $e = 1 <  j$, so $|(E_{r,s}^{-1})_{e+1,j+1}| = O(q^{-(j-e)s})$. Thus, we have
\[B_{n,r,j,e} =  O \lt( \sum_{1 \leq s < n/j} q^{-(j-1)s}  \frac{(\log (n-js))^{r-j-1}}{n-js}\rt).\]

Hence, regardless of the value of $e$, we have
\begin{align*}
B_{n,r,j,e} &=  O \lt( \sum_{1 \leq s < n/j} q^{-s}  \frac{(\log (n-js))^{r-j-1}}{n-js}\rt) \\
&=  O \lt( \sum_{1 \leq s < n/j} q^{-s}  \frac{(\log (n-js))^{r-j-1}}{n-js}\rt) \\
&=  O \lt( \sum_{1 \leq s < n/j} q^{-s}  (\log (n-js))^{r-j-1}\rt) \\
&=  O (\log (n-js))^{r-j-1}) \\
&=  O (\log (n-js))^{r-3}),
\end{align*}

because $j \geq 2$. This finishes the proof.
\end{proof}

\

\section*{Acknowledgments}

\hspace{3mm} We thank Yifeng Huang, Nathan Kaplan, Ofir Gorodetsky, and Michael Zieve for helpful conversations. G. Cheong was supported by NSF grant DMS-1162181 and the Korea Institute for Advanced Study for his visits to the institution regarding this research. We deeply appreciate the referee for quick and detailed comments on the previous draft of this paper. J. Lee was supported by a KIAS Individual Grant (MG079602) at Korea Institute for Advanced Study. H. Nam was supported by the National Research Foundation of Korea (NRF) grant funded by the Korea government (MSIT) (No. 2021R1F1A106231911). M. Yu was supported by a KIAS Individual Grant (SP075201) via the Center for Mathematical Challenges at Korea Institute for Advanced Study. He was also supported by the National Research Foundation of Korea (NRF) grant funded by the Korea government (MSIT) (No. 2020R1C1C1A01007604). This research was also supported by the Yonsei University Research Fund of 2022-22-0125.


\end{document}